\documentclass[12pt,a4paper]{article}
\usepackage{amsmath,amsthm,amsfonts,amssymb}
\usepackage{verbatim}
\usepackage[latin1]{inputenc}
\theoremstyle{plain} 
\newtheorem{thm}{Theorem}[section]

\newtheorem{lem}[thm]{Lemma}

\theoremstyle{definition}
\newtheorem{defn}{Definition}[section]

\theoremstyle{remark}
\newtheorem{rem}{Remark}[section]

\newcommand{\thmref}[1]{Theorem~\ref{#1}\xspace}
\newcommand{\secref}[1]{\S\ref{#1}\xspace}
\newcommand{\lemref}[1]{Lemma~\ref{#1}\xspace}

\newif\ifpdf
\ifx\pdfoutput\undefined
        \pdffalse 
\else
        \pdfoutput=1 
        \pdftrue
\fi

\ifpdf
        \usepackage[pdftex]{graphicx}
\else
        \usepackage[dvips]{graphicx}
\fi
 
\ifpdf
        \DeclareGraphicsExtensions{.pdf, .ppm, .jpg, .tif}
\else
        \DeclareGraphicsExtensions{.eps, .jpg}
\fi

\newcommand{\mapping}[1]{\xrightarrow{#1}}

\ifx\undefined\xspace \usepackage{xspace} \fi

\newcommand{\qr}[1]{\eqref{#1}} 
\newcommand{\st}{\mathrel{\,:\,}}
\newcommand{\ie}{\textit{i.e.}\xspace} 
\newcommand{\eg}{\textit{e.g.}\xspace}

 \newcommand{\RR}{{\mathbb R}}

\newcommand{\gparens}[3]{\left#1 #2 \right#3} 
\newcommand{\parens}[1]{\gparens({#1})} 
\newcommand{\brackets}[1]{\gparens\{{#1}\}} 
\newcommand{\hakparens}[1]{\gparens[{#1}]}

\newcommand*{\Setof}[1]{\,\brackets{#1}\,}

\renewcommand{\Pr}[2][]{\probopi{Pr}{#1}{#2}} 
 
\newcommand{\Ex}[2][]{\probopi{E}{#1}{#2}}

\newcommand{\Ordo}[1]{{O\parens{#1}}} 
\newcommand{\ordo}[1]{{o\parens{#1}}} 
\newcommand{\OrdoOmega}[1]{{\varOmega\parens{#1}}} 
\newcommand{\OrdoTheta}[1]{{\Theta\parens{#1}}} 
\newcommand{\ordoomega}[1]{{\omega\parens{#1}}}

\renewcommand{\L}{\mathsf{L}}
\newcommand{\J}{\mathsf{J}}

\newcommand{\e}{\varepsilon}

\newcommand{\Cols}{\mathcal{K}}
\newcommand{\Diags}{\mathcal{D}}

\newcommand{\Rows}{\mathcal{R}}
\newcommand{\Syms}{\mathcal{S}}
\newcommand{\X}{X}

\newcommand{\g}{\gamma}

\newcommand{\triplet}[3]{{#1,#2,#3}}

\newcommand{\ikg}[1][i]{\triplet{#1}{k}{\g}}

\newcommand{\RKS}{\mathcal{X}}

\newcommand{\UI}{{[0,1]}}

\newcommand{\p}{\mathsf{p}}
\newcommand{\q}{\mathsf{q}}

\newcommand{\ii}[1]{^{#1}}

\newcommand{\tr}[1][t]{\ii{#1+1}}
\renewcommand{\t}[1][t]{\ii{#1}}

\newcommand{\row}{i}

\newcommand{\K}{\mathsf{K}}

\newcommand{\CF}{\mathcal{F}}

\renewcommand{\Pr}[1]{{\mathbb{P}}\!\brackets{#1}}
\renewcommand{\Ex}[1]{{\mathbb{E}}\!\hakparens{#1}}
\newcommand{\Prt}[2][t]{{\mathbb{P}_{#1}}\!\brackets{#2}}
\newcommand{\Ext}[2][t]{{\mathbb{E}_{#1}}\!\hakparens{#2}}

\newcommand{\Good}{\Gamma}
\newcommand{\Lats}{\mathfrak{L}}

\newcommand{\bool}{{\{0,1\}}}

\newcommand{\T}{T}
\newcommand{\ctime}{\tau_{c}}

\newcommand{\Time}{{[0,m]}} 

\newcommand{\coord}[1][]{\iota_{#1}}
\newcommand{\RS}{{\Rows\Syms}}
\newcommand{\RK}{{\Rows\Cols}}
\newcommand{\KS}{{\Cols\Syms}}
\newcommand{\DS}{{\Diags\Syms}}

\newcommand{\M}{\mathsf{M}}

\newcommand{\ro}[1][]{\ell_{#1}}
\newcommand{\rok}{\ro[\KS]}
\newcommand{\rod}{\ro[\DS]}

\newcommand{\tp}{\p_\ell} 
\newcommand{\Tc}[1][\ell]{\tau_{#1}}

\newcommand{\asymabbr}[1]{\operatorname{\mathfrak{#1}}}
\newcommand{\pmax}{\asymabbr{p}}
\renewcommand{\aa}{\asymabbr{a}} 
\newcommand{\qpd}{\asymabbr{b}}

\begin{document}

\title{Orthogonal latin rectangles}
\date{Sep 21, 2004}

\author{
\large Roland H{\"a}ggkvist \\
\small Matematiska institutionen, Umeå universitet,\\ 
\small S-901 87 Umeå, Sweden \\
\normalsize and\\
\large Anders Johansson \\
\small N-institutionen, Högskolan i Gävle\\ 
\small S-801 76 Gävle, Sweden\\
\small Email: {\tt ajj@hig.se, rolandh@math.umu.se}
}

\maketitle
\begin{abstract}
  We use a greedy probabilistic method to prove that for every
  $\epsilon > 0$, every $m\times n$ Latin rectangle on $n$ symbols has
  an orthogonal mate, where $m=(1-\epsilon)n$.  That is, we show the
  existence of a second latin rectangle such that no pair of the $mn$
  cells receives the same pair of symbols in the two rectangles.
\end{abstract}    

\section{Introduction}
        
This paper was inspired by a problem posed by Anthony J. W. Hilton at
the thirteenth British Combinatorial Conference 1991 \cite{Hilton}.
The problem is:
\begin{quote}
  Let $R$ be an $n\times 2n$ Latin rectangle on $2n$ symbols. A
  partial transversal $T$ of size $s$ of $R$ is a collection of $s$
  cells, no two in the same row or column, and no two containing the
  same symbol. Is it true that $R$ can be expressed as the union of
  $2n$ partial transversals of size $n$?

  An equivalent formulation: Call two $n\times 2n$ Latin rectangles
  $R,S$ on the same set of symbols \emph{orthogonal} if the pairs
  $(r_{ij},s_{ij})$, for $i=1,\dots,n$ and $j=1,\dots,2n$, are all
  distinct. Does every $n\times 2n$ Latin rectangle have an orthogonal
  mate?
\end{quote}
The updated problem list from the British Combinatorial Conferences
from number twelve and upwards can be found electronically as a link
on the homepage of the British Combinatorial Conference.

While it is quite easy to see that every $n\times 4n$ Latin rectangle
has an orthogonal mate we know of no argument that solves the problem
for $n\times 3n$, when $n=100$, say. No doubt such an argument can be
found eventually.
 
In the current paper we use probabilistic methods to prove a far
stronger statement than Hilton proposed, but only valid for large $n$,
namely that for every $\e > 0$, every $(n-\e n) \times n$ Latin
rectangle has an orthogonal mate for large enough $n$.

We know of no example of an $(n-1)\times n$ Latin rectangle without an
orthogonal mate, but would not be too surprised if such an example
could be constructed. This ties up with well-known conjectures and
results concerning the length of partial transversals in latin
squares.  Recall that Ryser \cite{Koksma}, Brualdi \cite[s. 103]{DK}
and Stein \cite{Stein} have conjectures (in particular Stein has much
stronger conjectures, one of which was refuted by Drisko
\cite{Drisko}) which imply that every $(n-1)\times n$ latin rectangle
has a transversal of length $n-1$. In this context we also recall some
standard results on the length of partial transversals in latin
squares, to viz: every $n\times n$ latin square has a partial
transversal of length at least $n-\sqrt n$ (proved by Woolbright
\cite{Woolbright}, and Brouwer, de Vries and Wieringa \cite
{Brouwer}), and $n-5.53(\log n)^2$ proved by Shor \cite{ Shor}.

\subsection{The result once again}
We consider \emph{$m\times n$-latin rectangles} on $n$ symbols, $n$
columns and $m$ rows, \ie, an assignment to any cell in an $m\times
n$-table one of $n$ symbols such that each symbol occur exactly once
in each row and at most once in each column.

Two latin $m\times n$-rectangles, $\L$ and $\J$, are \emph{orthogonal}
if the following holds: For any two colours $\alpha,\beta$ the
colour-classes $\L^{-1}(\alpha)$ and $\J^{-1}(\beta)$ intersect in at
most one element.  Equivalently, each colour class $\L^{-1}(\alpha)$ is
a \emph{transversal} of $\J$ and vice versa.
        
\begin{thm}\label{main1} For every $\e > 0$ there is an $n_0 =
  n_0(\e)$ such that to any $m\times n$-latin rectangle $\J$, $n \geq
  n_0$ and $m = n (1-\e)$, there is an orthogonal companion $\L$.
\end{thm}
\begin{rem}
With some extra effort it is perhaps possible to prove
\thmref{main1} for all $\e = \ordoomega{n^{-1/2}}$.  However, it will
be clear from the proof that
in order to reach $\e \leq n^{-1/2}$ some new ideas must be
found, if indeed the theorem is valid in this range. 
\end{rem}

The basic method is related to nibble-methods used to colour graphs
having ``near disjoint'' cliques.  An orthogonal companion $\L$ of
$\J$ can be thought of as an $n$-colouring of the graph where the $mn$
cells are vertices and where each row, column and $\J$-colourclass
make up a clique, \ie a complete induced subgraph.  The recent
monograph of Reed and Molloy \cite{MolloyReed} contains many result in
this area, \eg J. Kahn's result \cite{kahn:asym} on edge-colourings of
near-disjoint hypergraphs.

There are elements of this proof that we belive are new. First of all,
we have a distinguished parallel class of cliques of size $n$ ---
corresponding to the rows --- that we colour with $n$ colours.  In
other words we have no slack here.  

In the proof, we construct an orthogonal companion in a random greedy
manner by adding one row at a time. We use a process, $\q^t$,
$t\in\Time$, of ``fractional latin rows'' to guide the greedy
extensions, so that
$\q^t(t+1,\cdot,\cdot)\in\RR^{\Cols\times\Syms}$ gives the
expectation of the row, $t+1$, added at time $t$. We maintain the
``legality'' of $\q^t$ by setting $\q^t(\ikg)=0$ if the cell $(i,k)$
belongs to a $\J$-colourclass or column already coloured with symbol
$\g$.

By analysing the time evolution for certain statistics of the random
process $\q\t$, we deduce that, with positive probability, $\q\t$ can
be legally maintained for all times $t=0,1,\dots,m-1$ so that the
latin rectangle $\L$ is constructed at time $m$.

\subsection{Rows, columns, cells, diagonals and points}

We will think of a vector $f$ in Cartesian space $\RR^A$ as a
real-valued mapping $f$ from the index set $A$. Pointwise relations
extends to relations between vectors in the natural way, \eg $f\leq g$
means that $f(a)\leq g(a)$ for all $a\in A$.

Let $\Rows = [1,m] := \{1,2,\dots,m\}$ denote the set of \emph{rows},
$\Cols$ denote the set of \emph{columns} and $\Syms$ the set of
\emph{symbols}.  Thus, $|\Rows| = m = n(1-\e)$ and $|\Cols| = |\Syms|
= n$.  We refer to elements of $\Rows\times\Cols$ as \emph{cells}.
Let $\J$ be the given latin rectangle from \thmref{main1}.  A
\emph{diagonal} is a set of cells assigned a common colour by $\J$ and
the family of diagonals is denoted by $\Diags$. The elements of $\RKS
:= \Rows\times\Cols\times\Syms$ we refer to as \emph{points}. If
nothing else is stated, we assume that the variable $i$ refer to a
row, the variables $k,l$ refer to columns and a variable $\gamma$
refers to a symbol. The variables $x$ and $y$ will be preferred for
points.


We have, for $\mathcal{P} = \Rows,\Cols,\Diags,\Syms$, mappings $\RKS
\mapping{\coord[\mathcal{P}]} \mathcal{P}$ assigning to each point the
unique row, column, diagonal or symbol to which it belongs.  We
usually say that a point $x\in\RKS$ \emph{belongs to} the corresponding
row, column, diagonal or symbol $\alpha \in \mathcal{P}$, when
$\coord[\mathcal{P}](x) = \alpha$.

A \emph{line} is a set of points with two of these coordinates fixed,
\ie a line is the set of the form
$(\coord[\mathcal{P}]\times\coord[\mathcal{Q}])^{-1}(\alpha,\beta)$
with $\alpha \in \mathcal{P}$ and $\beta \in \mathcal{Q}$.  We
introduce for any pair of two distinct coordinates $\mathcal{P}$ and
$\mathcal{Q}$ the mapping $\ell_{\mathcal{PQ}}$ assigning to each
point the corresponding line to which it belongs.  More precisely, we
let
\begin{equation}\label{elldef}
  \ell_{\mathcal{PQ}}(x) := 
  (\coord[\mathcal{P}]\times\coord[\mathcal{Q}])^{-1}
  (\coord[\mathcal{P}]\times\coord[\mathcal{Q}](x)).
\end{equation}
The collection of lines make up a linear hypergraph on the points
$\RKS$, \ie for every pair of distinct points $x,y\in \RKS$ there is
at most one line containing them.

\subsubsection{Latin rectangles}

A \emph{rectangle} can now be identified with a $\bool$-valued vector
$\L \in \bool^\RKS$ in the obvious way: For a point $x=(\ikg)$,
$\L(x)=1$ if $\g$ is the symbol assigned to cell $(i,k)$ and $\L(x) =
0$ otherwise.

A rectangle $\L \in \bool^{\RKS}$ is a \emph{latin rectangle
  orthogonal to $\J$} exactly when the relations
\begin{align*}\tag{L}
 \sum_{y \in \ell_\RK(x)} \L(y) &=  1, &
  \sum_{y \in \ell_\RS(x)} \L(y) &=  1, \label{local}\\[2ex]
            \tag{C}
            \sum_{y\in \ell_\KS(x)}  \L(y) &\leq 1, & 
            \sum_{y\in \ell_\DS(x)}  \L(y) &\leq 1, \label{central}
\end{align*}
hold for all $x \in \RKS$.  The relations in \qr{local} and
\qr{central} define a polytope $\Lats \subset \UI^{\RKS}$ so that a
latin rectangle $\L$ is a $\bool$-valued element of this polytope.
For our purposes, \emph{rational latin rectangles} orthogonal to $\J$
are vectors in $\Lats$.  The constraints in \qr{local} are
\emph{local} to each row since they concern lines contained in rows.
The constraints in \qr{central} are then \emph{central} constraints
since they concern lines transversal to the rows.


\subsection{The greedy latin rectangle process}

We now give a birds eye view of the proof. The probabilistic
terminology used regarding vector-valued random processes is made
precise in section \secref{probterm} below.  Our purpose is to
construct an increasing random process, a \emph{greedy rectangle
  process}, $\L\t\in\bool^\RKS$ of partial $\J$-orthogonal latin
rectangles that proceed row-wise: 
Initially, $\L^0 \equiv 0$ and at each tick of the clock, \ie when $t
\mapsto t+1$, we extend --- if the situation allows it --- the partial
latin rectangle $\L\t$ to a partial latin rectangle $\L\tr$ having the
row $t+1$ added to the latin rectangle.  The time variable $t\in\Time
:= \Setof{0,\ldots, m}$ thus corresponds to rows being added to the
rectangle.  The process is \emph{successful} if $\L^{m}$ actually
produce a full $\J$-orthogonal latin rectangle.

It is quite easy to see that such a greedy rectangle process should
always be successful as long as $m\leq n/4$.  To see this, note that
the legal choices of the added row $t+1$ are, by the local
constraints, given by matchings in a ``legality graph'', which is the
balanced bipartite graph consisting of those symbol-column pairs for
row $t+1$ that are not in conflict with any previously added row due
to central constraints. Moreover, each previously added row can
exclude at most two symbols for the column $k$ on row $t+1$; one
symbol on the same column and one symbol on the same diagonal.
Similarily, each added row excludes at most two possible columns for
any symbol $\gamma$.  Thus, since $t < m \leq n/4$ are previously
added, the legality graph will have minimum degree at least $n-2t \geq
n/2$ and a well known degree condition based on Halls theorem ensures
the existence of a legal matching for row number $t+1$.

This naive argument can be extended significantly when the obtained
legality graph is sufficiently random-like to ensure the existence of
a perfect matching for degrees well-below $n/2$. To achieve this, we
need to introduce some probabilistic tools. 

A central idea is to let the greedy rectangle process $\L\t$ be
``guided'' by a Markov process $\p\t \in \UI^\RKS$, $t\in\Time$. We
refer to $\p\t \in \UI^\RKS$ as a \emph{state}.  The initial state is
the uniform vector $\p^0 \equiv \frac 1n$. The relationship between
the processes $\L\t$ and $\p\t$ is that, at time $t$, $\p\t(x)$
approximately gives the expectation of $\L\t[s](x)$, for points $x$
belonging to rows that are coloured at time $s>t$.

Care must therefore be taken in the construction of $\p\t$, so that
$\L\t$ never violates the local and central constraints, \qr{local}
and \qr{central}.  We defer the exact definition $\p\t$ to section
\secref{pconstr} below. We note here that the construction of $\p\t$
ensures that the central constraints \qr{central} are never violated
by $\L\t$: If a cell $(t+1,k)$ in the active row is assigned the
colour $\gamma$ at time $t$ then we remove the possiblity that any
cell in the same column or diagonal later gets colour $\gamma$. Hence,
we must ``kill'' all points $y$ belonging to the central lines going
through the point $(t+1,k,\gamma)$, that is, we set
\begin{equation}\label{zzkill}
\p\tr(y)=\p\t[t+2](y)=\dots=0,
\end{equation}
for all $y$ belonging to such a central line.

Given $t\in\Time$, we define a region $\Good \subset
\UI^\RKS$, where $\p\t \in \Good$ should be interpreted as stating
that $\p\t$ is a ``good state''.  The exact definition of $\Good$ is
deferred to \secref{Gammadef} below, but we mention that $\Good$ is
defined by three collections of inequalities: The first group of
inequalites, \qr{pbnd1}, bounds the size of the individual values
$\p\t(x)$ while the second group,\qr{norm1}, states that $\p\t$ almost
should satisfy the local constraints \qr{local}. 

Note that, for a fixed row $i \in\Rows$, the local constraints given
by \qr{local}, defines a polytope $\Lats_i$ in
$\RR_+^{\coord[\Rows]^{-1}(i)} \cong \RR_+^{\Cols\times\Syms}$ which can
be interpreted as the polytope of \emph{rational (perfect) matchings}
in the complete bipartite graph $K(\Cols,\Syms)$ and we refer to
$\bool$-valued vectors in $\Lats_i$ as \emph{matchings} on that
row. 
\begin{lem}\label{goodadm}
  If $\p\t \in \Good$ then, for each row $\row\in \Rows$, there is a
  rational matching $\q\t_\row \in \Lats_\row$ such that for all
  $(k,\g)\in\Cols\times\Syms$
\begin{equation}\label{qpeta}
  \q^t_{\row}(\row,k,\g) \leq (1+\qpd)\p\t(\row,k,\g),
\end{equation}
where $\qpd$ is an abbreviation of the asymptotic expression
$\Ordo{\sqrt{n^{-1/2}\,\log n}}$. {\rm (}\/See \qr{asymabbr}
below.\/{\rm )}
\end{lem}
We prove this Lemma in section \secref{proofadmissible} using the
Ford-Fulkerson Theorem; in the proof, the third and last group of
inequalities, \qr{qrnd1}, which give a ``quasi-random'' property of
$\p\t$, are central for the construction.

Now recall the well-known characterization by Birkhoff
\cite{birkhoff} stating that any rational matching
$\q_\row\in\Lats_\row$ can be expressed as a convex combination
$q_\row = \sum_{M} c_M M$ of matchings $M\in\Lats_\row$. 
By interpreting the convex coefficients $c_M$ as probabilities, where
we pick the matching $M$ with probability $c_M$, Birkhoffs theorem can
also be given the following formulation: Given any rational matching
$\q_\row\in\Lats_\row$, it is always possible to find a \emph{random} matching
$\L_\row\in\Lats_\row$, such that the \emph{expectation} of $\L_\row$
equals $\q_\row$. 

Therefore, modulo the precise definition of $\Good$
and $\p\t$, we can define the greedy latin rectangle process $\L\t$ by
iterating the following procedure for $t=0, 1, \dots, m-1$.
\begin{description}
\item[Extend] If $\p\t\in\Good$ then choose a rational matching
  $\q_{t+1}\t$ on row $t+1$ which satisfies \qr{qpeta}. Then
  draw, using the random construction implied by Birkhoffs theorem, an
  extension $\L\tr$ such that
  \begin{equation}\label{zzequation}
    \Ext{\L_{t+1}\tr} = \q\t_{t+1}.
  \end{equation}

  The new state $\p\tr$ is then constructed from $\p\t$, $q_{t+1}$
  and $\L\tr$ according to the construction in \qr{pdef} below.

\item[Stop] If $\p\t \not\in\Good$ then simply let $\L\t[s] = \L\t$
  and $\p\t[s] = \p\t$ for all $s$, $t<s\leq m$. The greedy latin
  rectangle is then said to be unsuccessful.
\end{description}
On account of the killing mechanism \qr{zzkill}, the bound \qr{qpeta}
and the property \qr{zzequation}, the construction ensures that the
central constraints are never violated by $\L\t$.  Thus, if $\p\t$
stays in $\Good$, the process produces an orthogonal companion
$\L^{m}$ to $\J$ at time $m$ and the probability of an unsuccessful
rectangle process is the probability that $\p\t$ leaves $\Good$ for
some $t\in\Time$.  The proof is thus concluded, if we, after properly
describing the construction of $\p\t$ and the definition of $\Good$
and proving \lemref{goodadm}, in addition, prove the following lemma.
\begin{lem}\label{stable} For all $t=1,2,\dots,m-1$ we have that
\begin{equation}\label{zzzuuu}
	\Pr{\p\t \in \Good} = 1 - n^{-\ordoomega 1}. 
\end{equation}
\end{lem}
Note that \qr{zzzuuu} implies that the probability that $\p\t$ stays
inside $\Good$ at all times $t$ is of probability of order
$1-n^{-\ordoomega 1+1}=1 - n^{-\ordoomega 1}$. 

\subsection{Probabilistic preliminaries and asymptotic notation}
\label{probterm}

\subsubsection{Asymptotic notation}

We will use the standard asymptotic notation, $\Ordo{\cdot}$,
$\ordo{\cdot}$, $\ordoomega{\cdot}$, $\OrdoOmega\cdot$, etc., where
all are interpreted as asymptotic estimates relative the limit
$n\to\infty$. That is, 
$f = \Ordo g$ if and only if $\limsup_{n\to\infty} |f/g| < \infty$,
$f = \ordoomega g$ if and only if $\liminf_{n\to\infty} |f/g| = \infty$,
$f = \OrdoOmega g$ if and only if $\limsup_{n\to\infty} |g/f| < \infty$,
$f = \ordo g$ if and only if $\limsup_{n\to\infty} |f/g| = 0$ and 
$f = \OrdoTheta g$ if and only if $\limsup_{n\to\infty} (|f/g| + |g/f|)
<\infty$. 

Such asymptotic expressions are often used to estimate
the components of vectors and, of course, if we have such a
local quantity expressed in terms of some asymptotic expression, then
all implicit constants in the asymptotic expression are assumed to be
independent of the particular point, row, column et cetera at which the
local quantity is defined.

Since some asymptotic expressions are extensively used, we also
introduce the following {\em abbreviations} of asymptotic expressions
\begin{equation}\label{asymabbr}
        \begin{aligned} 
        \aa &:= \Ordo{ \frac {\log n}{\sqrt{n}} }, \quad 
        \qpd &:= \Ordo{n^{-1/4}(\log{n})^{1/2}}, \quad
        \pmax &:= \Ordo{ \frac {1}{n} }.
        \end{aligned}
\end{equation}
Note that $\sqrt{\aa} = \qpd$.

\subsubsection{Probabilistic terminology and notation}

The proof will use a dynamic probabilistic method, so we introduce
some concepts and terms from probability theory.  We will construct a
\emph{filtered probability space} $(\Omega,\Pr{\cdot},\CF^t)$ with a
discrete time-variable $t$ taking values in $\Time :=
\Setof{0,1,\dots,m}$ and we can assume that the probability space we
work with is finite.  The finite algebras $\{\CF^t\}$, $t\in\Time$, is
an increasing sequence of subsets of $2^\Omega$ with $\CF^0 =
\Setof{\emptyset,\Omega}$.  In our case, $\CF^t$ captures the random
operations used for adding the first $t$ rows to a latin rectangle. A
random variable is \emph{determined at time $t$} if it is
$\CF^t$-measurable.

We will work with vector valued random variables without explicitly
noting this: A (vector-valued) random variable is a mapping $X: \Omega
\mapping{} \RR^A$ from $\Omega$ to a Cartesian space $\RR^A$. For our
purposes, a \emph{stochastic process} is a function $X: \Omega\times
\Time \to \RR^{A}$. We write $X(\omega,t)$ simply as $\X^t$. Functions
defined for $t\in \Time$ will generally have the variable $t$ as a
superscript.  We suppress the dependence on $\omega \in \Omega$
for the random entities.

A process $X^t$ is \emph{adapted} if the value of $X^t$ is determined
at time $t$ and we usually assume this to be the case without further
notice.  A process $\X\t$ is \emph{previsible} if the value of $X^t$
is determined at time $t-1$.

The expectation operator $\Ex{\cdot}$ and the probability $\Pr{\cdot}$
refers to the unconditional probability, while the temporal
conditional expectation is denoted by $\Ext{f} = \Ex{f \mid \CF^t}$
and the conditional probability by $\Prt{A} = \Pr{A \mid \CF^t}$.

The various expectation operators apply to \emph{vectors} so that if
$F$ is a random vector taking values in a Cartesian space $\RR^{A}$ we
mean by $\Ex{F}$ the vector in $\RR^A$ given by $\Ex{F}(a) =
\Ex{F(a)}$, $a\in A$.

An adapted process $X\t$ is a \emph{martingale (super-martingale or
  submartingale)} if $\Ext{X\tr} = X\t$ ($\Ext{X\tr} \leq X\t$ or
$\Ext{X\tr} \geq X\t$).

A \emph{stopping time} is a random time $\tau: \Omega \to \Time \cup
\Setof{\infty}$ such that the event $\{\tau \leq t\} \in \CF^t$.  We
will work with \emph{vectors of random times}, \ie, mappings $\tau:
\Omega \times A \to \Time \cup \Setof{\infty}$.

Given a vector-valued process $X\t \in \RR^{A}$, $t\in\Time$, then
$X^{t\wedge \tau}$ is the process whose value at $a\in A$ at time $t$
is $X^t(a)$ if $t \leq \tau(a)$ and $X^{\tau(a)}(a)$ otherwise. (We
use $s\wedge t$ as a shorthand for $\min\{s,t\}$.)

If $X\t$ is adapted and $\tau$ is a vector of stopping times then
$X^{t\wedge \tau}$ is adapted.  We say that an adapted process $X\t$
is \emph{stopped} at a vector of stopping times $\tau$ if $X\t =
X^{t\wedge \tau}$. If $X\t$ is a supermartingale and $\tau$ is a
vector of stopping times then the process $X^{t\wedge\tau}$ is a
supermartingale.

Before proving \lemref{goodadm} and \lemref{stable} in the following
sections, we first proceed to define the process $\p\t$ rigorously as
well as the good set $\Good$.

\subsection{The construction of $\p\t$}\label{pconstr}

For ease of notation, we define a deterministic vector $\ctime: \RKS
\to\Time$ of \emph{colouring times}, so that, for a point $x =
(\ikg)$, $\ctime(x) = i-1$ is the time when the row $i$ is added to the
latin rectangle process.

\subsubsection{The killing}
Consider a point $x = (\ikg) \in \RKS$ such that $\ctime(x) \not= t$.  For
a central line $\ell = \ell_\KS$ or $\ell= \ell_\DS$, let $\ro^t(x)$
be the unique point on the active row $t+1$ lying in the line
$\ell(x)$, $x \in \RKS$. If any of $\L\tr \circ \rok^t(x)$ and $\L\tr
\circ \rod^t(x)$ take the value $1$ then the point $x$ is
\emph{killed} at time $t$.

The two ``projections'' of $x$, $\rok^t(x)$ and $\rod^t(x)$, are
distinct points and belong to a common local line
$\ell_{\RS}(\rok^t(x)) = \ell_{\RS}(\rod^t(x))$ in the row $t+1$.
The indicators $\L\tr \circ \rok^t(x)$ and $\L\tr \circ \rod^t(x)$ can
therefore never take the value $1$ simultaneously and we can write
\begin{equation}\label{Kdef}
    \K\tr(x) = \L\tr \circ \rok^t(x) + \L\tr \circ \rod^t(x),
\end{equation}
for the indicator $\K\tr(x) \in \bool$ of the event that the point $x$
is killed.

\subsubsection{The construction of $\p\t$}
Initially, we set $\p^0(x) \equiv \frac 1n$ for all $x\in\RKS$.  We
define the global stopping time $\T$ marking exit from $\Good$, \ie
\begin{displaymath}
  \T := \min \Setof{t \in [0,m-1] \st \p^t \not\in \Good} \cup\{\infty\}.
\end{displaymath}
For $t > \T$ the greedy process is thus in effect ``stopped'' and the
greedy random colouring has failed if $T < \infty$.

Therefore, for $t < \T$ and $x\in\RKS$ such that $\ctime(x) > t$,
define
\begin{equation}
  \p\tr(x)  := \p\t(x) \cdot \frac{1-\K\tr(x)}{\Ext{1-\K\tr(x)}}
        \label{pdef}
\end{equation}
and for $t \leq \T$ and $t \geq \ctime(x)$, let $\p\tr(x) = \p\t(x)$.

By the definition of $\p\tr$ above, the process $\p\t$ is
\emph{stopped} both at the global stopping time $\T$ and at the
deterministic stopping time vector $\ctime: \RKS \to \Time$, \ie,
\begin{math}
  \p\t = \p^{t\wedge\ctime\wedge T}.
\end{math}

\subsubsection{Some properties of $\p\t$}
For $t<\T$, \ie as long as $\p\t\in\Good$, we can, by Definition
\ref{defgood} below, assume that
\begin{equation}\label{pleqpmax}
	\p\t \leq \pmax.
\end{equation}
with the asymptotic abbreviation $\pmax=\Ordo{1/n}$ as in
\qr{asymabbr}. (We understand that a relation like \qr{pleqpmax} holds
with the same implicit constant at all points $x\in\RKS$. )

Since the vector $1-\K\tr\in\{0,1\}^\RKS$ indicates survival, the
definition ensures that all points $x\in\RKS$ such that $\ctime(x) >
t+1$ and $\p\tr(x) > 0$ can be used to extend $\L\tr$. It follows that
$\L\t$ is indeed a process of legal (partial) $\J$-orthogonal latin
rectangles.

Note that, from \qr{qpeta}, \qr{Kdef} and \qr{pleqpmax}, we have that
\begin{equation}\label{exstr}
    \begin{aligned}
      \Ext{1-\K\tr} &= 1 - \q\circ\rok\t - \q\circ\rod\t\\
      &= 1 - \p\circ\rok\t - \p\circ\rod\t - \qpd \pmax \\
      &= (1 + \pmax)^{-1}.
    \end{aligned}
\end{equation}
Hence
\begin{equation}\label{pincr}
    \p\t(x) \leq \p\tr(x) \leq \p\t(x) \cdot (1+\pmax), 
\end{equation}
unless $x$ is \emph{killed}, \ie unless $\p\tr(x) = 0$.

By the construction \qr{pdef}, the process $\p\t$ is a
\emph{martingale}, so that
\begin{equation}\label{martingale}
        \Ext{\p\tr} = \p\t.
\end{equation}
Relations \qr{pincr} and \qr{martingale} will be essential in deriving
the concentration results upon which the proof is founded.

\subsection{The definition of $\Good$}\label{Gammadef}

The process $\p\t$ is controlled by keeping a set of local
inequalities alive through the iterations.  These local relations make
up the notion of ``goodness'' that we rely on throughout the
arguments.  First, for $x=(\ikg) \in \RKS$, let the inequality
\qr{pbnd1} be defined by
\begin{equation*}
        \tag{$A_{x}$} \p^t(x) \leq 1.1 \, \e^{-2}\, \frac 1n.  \label{pbnd1} 
\end{equation*}
Secondly, for a local line $\ell = \ell_\RK, \ell_\RS$ and $x\in
\RKS$, we say that \qr{norm1} hold if
\begin{equation*}
  \tag{$B_{\ell(x)}$} 
  (1 - \frac{\log n}{\sqrt n}) \leq \sum_{y\in\ell(x)} \p\t(y)   
  \leq (1 + \frac{\log n}{\sqrt n}).      \label{norm1}
\end{equation*}
Finally, for $k\not=l \in \Cols$ and $i \in \Rows$ let \qr{qrnd1} be
the ``quasi-random'' inequality
\begin{equation*}
  \tag{$C_{i,kl}$}    
  \sum_\g \p\t(\triplet ik\g) \, \p\t(\triplet il\g) 
  \leq \frac 1n\, (1 + \frac{\log n}{\sqrt n}). 
            \label{qrnd1}
\end{equation*}

\begin{defn}\label{defgood} 
  We say that the state $\p\t \in \UI^{\RKS}$ is \emph{good} if
  \qr{pbnd1}, \thetag{$B_{\ell_\RK(x)}$} and
  \thetag{$B_{\ell_\RK(x)}$} hold at all $x\in\X$ and, in addition,
  \qr{qrnd1} hold for all $k,l\in\Cols$ and $i\in\Rows$.  We write
  $\Good$ for the region $\Good \subset \UI^\RKS$ of good states.
\end{defn}
It is trivial to check that the initial state, $\p^0 = 1/n$, is an
element of $\Good$.

\subsubsection{Asymptotic relaxations}
The precise formulations of \qr{pbnd1}, \qr{norm1} and \qr{qrnd1} are
needed to make $\Good$ well defined, but that precision is otherwise
not critical. What we actually will use in the computations below
are the following less precise asymptotic statements
\begin{gather}
  \tag{$A$} \p^t(x) =  \pmax, \label{pbnd} \\
  \tag{$B$} \sum_{\ell(x)} \p\t  = 1 \pm \aa, \label{norm} 
\intertext{and}
  \tag{$C$} \sum_\g \p\t(\triplet ik\g) \, \p\t(\triplet il\g) \leq
  (1+\aa)\,\frac 1n. \label{qrnd}
\end{gather}

It should be noted that the states $\p^t$, if they eventually leave $\Good$ at
the time $\T$,  will stay quite close to $\Good$.
This is due to the fact that we stop $\p\t$ at time $t=\T$, with the
previous state $\p^{\T-1}$ being a good state. As a consequence, we
can use the somewhat relaxed bounds \qr{pbnd}, \qr{norm} and \qr{qrnd}
in our arguments, \emph{without considering if we are conditioning on
  $t < \T$ or not}, since these bounds then hold for all times
$t$. In particular, we can assume the relaxed bounds when we
later show that $\p\t$ with high probability stays inside $\Good$.
(For definiteness, one may choose to replace the implicit constants in
\qr{pbnd}, \qr{norm} and \qr{qrnd} with explicit constants slightly
larger or smaller than those used in the definition of $\Good$.)

In order to see that $\p^\T$ is close to $\Good$ in this sense, note
that from \qr{pincr} it is clear that the left hand side of \qr{pbnd}
can only increase with a fraction $1+\pmax$ at a time.  Similarily,
the left hand side of \qr{qrnd}, can at most increase by a factor
$(1+\pmax)^2=1+\aa$ at time. Finally, from the fact that at most two
terms (see argument preceding \qr{uuux} below) in the sums in
\qr{norm} can be killed at a time, it also follows that the left hand
side of \qr{norm} can only change with a fraction $1\pm \p$ at a time.

\section{The proof of \lemref{goodadm}}\label{proofadmissible}


The reason behind introducing the set $\Good$ is \lemref{goodadm}
stated in the introduction. This lemma shows the existence of $\q\t$
and hence lets us define the extension $\L\tr$ of $\L\t$ as long as
$\p\t\in\Good$. We now proceed to prove this lemma.
\newcommand{\CS}{\Cols\times\Syms}

We only have to look at one fixed row $i\in \Rows$ at a time.  Given a
good state $\p\t \in\Good$, let $p \in \UI^{\CS}$ be given by
\begin{equation}\label{defprow}
  p(k,\g) := \p\t(i,k,\g)/\sum_l \p\t(\triplet il\g), 
\end{equation} 
so that $p = (1\pm\aa)\, \p_i\t$ by \qr{norm} and $p$ is normalized at
each symbol, \ie
\begin{equation}
  \sum_k p(k,\g) = 1. \label{gn}
\end{equation}
We assume that $\p\t$ satisfies the inequalities \qr{pbnd}, \qr{norm}
and \qr{qrnd}, which translates to the following set of inequalities
\begin{align}
p(k,\g) &= \pmax\label{pbndx},
\\
\sum_{\g\in\Syms} p(k,\g) &= 1+\aa\label{normx},
\\
\sum_{\g\in\Syms} p(k,\g)p(l,\g) &\leq (1+\aa)/n\label{qrndx},
\end{align} 
for all values of $k,l\in \Cols$ and $\g\in\Syms$.

We want to show that for some $\eta = \qpd$ there exists a rational
matching, \ie a vector $q \in \UI^{\CS}$ such that for all $k$ and
$\g$,  $\sum_{k'} q(k',\g) = \sum_{\g'} q(k,\g')=1$, 
that in addition satisfies
\begin{displaymath}
  q(k,\g) \leq (1+\eta) \cdot p(k,\g), \quad \forall (k,\g) \in \CS.
\end{displaymath}
Such a rational matching can also be defined as a flow on a directed
graph from a source $s$ connected to all vertices $k\in \Cols$ to a sink
$t$ connected to all vertices $\g\in\Syms$. The flow should take the
value $1$ on each edge $sk$, $k\in \Cols$, and each edge $\g t$,
$\g\in\Syms$. On the remaining edges, of the form $k\g$, we prescribe
the capacities $c_{k\g}=(1+\eta)p(k,\g)$.
The Ford-Fulkerson theorem says that it is enough to show that for all
pairs of nonempty sets $A \subsetneq \Syms$ and $B \subsetneq \Cols$
we have
\begin{equation}
  2n - |A| - |B| + 
  (1+\eta) \sum_{\substack{\g\in A \\ k\in B}} p(k,\g) \geq n, 
        \label{ff1}
\end{equation}
since the left hand side is the capacity for a cut in a flow
defining a rational matching where the capacity of edge $k\g$ equals
$(1+\eta) \cdot p(k,\g)$.

Given an arbitrary pair of subsets $A$ and $B$ as above, we shall
prove that, for some $\eta = \qpd$ (with the implicit
constant independent of $A$ and $B$), the left hand side in \qr{ff1}
is not strictly less than $n$ if we assume that $p$ satisfies the
inequalities \qr{pbndx}, \qr{normx} and \qr{qrndx}. 

Hence, assume that \qr{ff1} does not hold and proceed to derive a
contradiction. It must obviously then be the case that
$|A|+|B|>n$ and both $A$ and $B$ must be non-empty.  Let $a :=
|A|/n$, $b := |B|/n$. Then, $a,b>0$ and  
\begin{equation}\label{apb}
  a+b > 1. 
\end{equation}

For $\g\in \Syms$ and $S \subset \Cols$, let
\begin{displaymath}
  p(S,\g) := \sum_{k \in S} p(k,\g), 
\end{displaymath}
and define
\begin{equation}\label{xxyydef}
  x := b - \frac 1{|A|} \sum_{ \g\in A } p(B,\g)
  \qquad\text{and}\qquad
  y := b - \frac 1{n-|A|} \sum_{ \g\in \Syms\setminus A } p(B,\g).
\end{equation}
It is easy to see that, by dividing both sides in \qr{ff1} by $n$ and
substituting for $a$, $b$ and $x$, the assumption that \qr{ff1} does
not hold is equivalent to
\begin{equation}
  (1-a)(1-b) + ab\eta - a x (1+\eta) < 0. \label{ff2}
\end{equation}

In order to contradict \qr{ff2}, it is enough to show that
\begin{equation}\label{ax}
  ax \leq ab\qpd,
\end{equation}
since we can take $\eta$ equal to, say, two times the $\qpd$-function
on the right hand side of \qr{ax}.

We claim it follows from \qr{apb} and \qr{pbndx}--\qr{qrndx} that
\begin{equation}
  |a x + (1-a) y| \leq ab \cdot \aa, \label{lr}
\end{equation}
and that
\begin{equation}
  a x^2 + (1-a) y^2 \leq a b \cdot \aa. \label{sqr}
\end{equation}

Postponing the proofs for the two relations \qr{lr} and \qr{sqr} until
later, we proceed to show that they imply \qr{ax}.  We divide into
two cases depending on the value of $a$: If $a \leq 1/2$ then \qr{sqr} gives
\begin{displaymath}
  a x^2 \leq ab \aa =  a b^2 \aa,
\end{displaymath}
where we use that \qr{apb} implies that $b > 1/2$ so that $\aa = b
\aa$.  Multiplying both sides with $a$ and taking the square root
gives \qr{ax} since $\sqrt{\aa} = \qpd$.

In the case $a > 1/2$ then \qr{sqr} gives that 
\begin{displaymath}
  (1-a) y^2 \leq ab \aa 
\end{displaymath}
and since $1-a < b$ by \qr{apb} we can multiply the left hand side by
$1-a$ and the right hand side by $b$. This gives
\begin{displaymath}
  (1-a)^2 y^2 \leq a b^2 \aa = a^2 b^2 \aa, 
\end{displaymath}
where the last inequality is due to $a > 1/2$. Taking the square root
of this, substituting in \qr{lr} and using the triangle inequality we get
\begin{displaymath}
  a|x| \leq ab\aa + (1-a) |y| \leq ab \aa + ab\qpd = ab\qpd 
\end{displaymath}
and \qr{ax} is proved.  \qed

(It follows that we can take $\eta$ to be four times the square root
of the maximum $\aa$-function found in \qr{lr} and \qr{sqr}.)
        

What remain now is to show that \qr{lr} and \qr{sqr} follows from
\qr{apb} and the goodness assumptions \qr{pbndx}, \qr{normx} and
\qr{qrndx}.  Note that, by \qr{apb},
$$ ab\aa \geq (1-b)b\aa \geq (1/2)\,\min\{b,1-b\}\aa = \min\{b,1-b\}\aa. $$
Consequently, it is enough to show \qr{lr} and \qr{sqr} for the right
hand side $\beta\aa$, $\beta \in \Setof{b,1-b}$, instead of $ab\aa$.

Moreover, since for each $\g$, $p(B,\g) = 1 - p(\Cols\setminus B,\g)$,
the value $x$ and $y$ will both change only in sign if we interchange
$B$ with $\Cols\setminus B$ in \qr{xxyydef}. Thus, if we do not use the
assumption \qr{ff2} (and \eg \qr{apb}) about $B$, we may freely
interchange $B$ and $\Cols\setminus B$. This means that we only have
to consider the case $\beta = b$.

We first show that $|ax+(1-a)y|\leq b\aa$: We have that $|A| (b-x) +
(n-|A|)(b-y)$ equals $\sum_\g p(B,\g)$ which, by \qr{normx}, is of
order $|B|(1+\aa)$.  Dividing by $n$ gives
\begin{displaymath}
  a(b-x) + (1-a)(b-y) = b + b\aa \iff |ax + (1-a)y| = b\aa.
\end{displaymath}
 
In order to see that that $ax^2+(1-a)y^2\leq b\aa$, we let
$z_{kl}=\sum_{\g\in\Syms} p(k,\g)\cdot p(l,\g)$ denote the left hand side of
\qr{qrndx} above.  Note that,
\begin{equation}\label{p2rel}
  \sum_\g p(B,\g)^2 = \sum_{\substack{k,l \in B \\ k\not= l}} z_{kl} + 
  \sum_{\substack{\g \\ k \in B}} p^2(k,\g).
\end{equation}

The last sum in \qr{p2rel} above is of the order $|B|\pmax$ by
\qr{pbndx} and the first sum on the right hand side is less than
\begin{math}
  |B|^2 \, (1+\aa) / n
\end{math}
by \qr{qrndx}.  Furthermore, the Cauchy-Schwartz inequality implies
that the left hand side of \qr{p2rel} is greater than
\begin{multline*}
    |A| (\sum_{\g\in A} p(B,\g))^2   + 
    (n-|A|) (\sum_{\g\not\in A} p(B,\g))^2  = \\
    |A|(b-x)^2 + (n-|A|)(b-y)^2. 
\end{multline*}
Thus, we have
\begin{displaymath}
  |A|(b-x)^2 + (n-|A|)(b-y)^2 \leq  |B|^2 \, (1+\aa)/n + |B| 
  \pmax
\end{displaymath}
and dividing by $n$ and expanding the squares gives the relation
\begin{gather*}
  a b^2 + (1-a) b^2 - 2 b \,(a x + (1-a) y) + a x^2 + (1-a) y^2 \leq
  b^2 + b^2 \aa + b \aa \intertext{which is equivalent to} a x^2 +
  (1-a) y^2 \leq b^2 \aa + b \aa + 2b\,(a x + (1-a) y).
\end{gather*}
The sought statement follows if we use \qr{lr} to estimate the last
term. \qed
 
\section{The proof of \lemref{stable}}

For the parameter $n$ tending to infinity, we stipulate that an event
has {\em very high probability} if it holds with probability having
asymptotic order $1 - n^{-\ordoomega{1}}$.  An inequality of the form
$X\t\leq f(t)$ that, for all $t \in \Time$, holds with very high
probability is said to be \emph{stable}.  Note that, since we only
consider $n^{\Ordo 1}$ different times $t$, a stable inequality holds
with very high probability \emph{simultaneously} for all $t\in\Time$.

In order to conclude the proof of the lemma, it suffices to show
that, for every possible value of $x,i,k,l,\g$ and $\ell$
\begin{equation}
  \text{the inequalities \qr{pbnd1}, \qr{norm1} and \qr{qrnd1} 
    defining ${\Good}$ are all stable}.
        \label{overgoal}
\end{equation}
Since the definition of $\Good$ considers $\Ordo{n^4}$ such
inequalities, the probability that $\p\t\not\in\Good$ is then shown to
be of order $\Ordo{n^{-\ordoomega{1} + 4}} = n^{-\ordoomega1}$.

In the first subsection, we state and prove a more general
``concentration lemma'' \lemref{concentration} and then we prove, in
three separate subsections, the stability of \qr{norm1}, \qr{qrnd1}
and \qr{pbnd1} where we regularly invoke \lemref{concentration}.  Of
the three, proving the stability of \qr{pbnd1} involves the most
complex argument, but is should be noted that analysing the structure
of the linear hypergraph given by the lines is also an essential
component to derive the other two statements.

Again it should be noted, as in \secref{Gammadef}, that, since we stop
$\p\t$ at time $\T$, we can use the bounds \qr{pbnd}, \qr{norm} and
\qr{qrnd} in our arguments to derive the stability of \qr{norm1},
\qr{qrnd1} and \qr{pbnd1}, without considering if we are conditioning
on $t < \T$ or not. It should be clear that we at no point make the
assumption that the process $\p\t$ is unstopped. In particular, the
conditions and conclusions of \lemref{concentration} work for stopped
processes as well as ``live'' processes.

\subsection{A concentration result}\label{azhoeff}

\newcommand{\drift}{\beta}
\newcommand{\dev}{\delta}

The following lemma is a consequence of \emph{Azuma--Hoeffding's
  inequality}, see \eg \cite{alsp:probmeth}.  In order to make the
subsequent invocations transparent, we put the lemma in a suitable form
and make no attempt to derive the best possible result.
\begin{lem}\label{concentration}
  Let $\xi$ and $\alpha^0,\alpha^1,\dots,\alpha^{m-1}$ be positive
  numbers such that $(\xi\,\sqrt n + a)\log n = \ordo 1$, where
  $a:=\sum_{t=0}^{m-1} \alpha^t$.  Let $\X^t \geq 0$, $t \in \Time$,
  be a positive process such that
  \begin{align}
    |X\tr - \Ext{ \X\tr }| &\leq \xi \, \max\{X\t,X^0\}
    \label{devv}
\intertext{and}
    \Ext{ \X\tr } - \X\t &\leq \alpha^t \, \max\{X\t,X^0 \}.
    \label{driftt}
  \end{align}
   Then, for all $t\in \Time$,
  \begin{equation}\label{uub}
    \Pr{X^t \leq   (1+\Phi)\, X^0} = 1 - n^{-\ordoomega1},
  \end{equation}
  for any $\Phi = \OrdoOmega{(\xi\sqrt{n} + a) \log n}$.  Furthermore,
  if $\X\t$ is a \emph{martingale} (in which case $a = 0$) then the
  reverse inequality
  \begin{equation}\label{lub}
    X^t \geq  (1 - \Phi) \, X^0 
  \end{equation}
  is stable as well.
\end{lem}

\begin{proof}[Proof of \lemref{concentration}]
Note that the inequalities \qr{devv} and \qr{driftt} are still valid
  if we stop the process at any stopping time $\tau$.  
If we take
  $\tau$ to be the first time that $X^t > 2 \cdot X^0$ then, since
$$ \X^{s+1}  \leq (1 + \xi + \alpha^s)\,X^s=(1+ \ordo{1})\,X^s$$
we can assume that the stopped process $X^t=X^{t\wedge \tau} \leq 3
X^0$.  It is therefore enough to prove the stability of \qr{uub} with
the additional assumption that $X\t \leq 3 \cdot \X^0$ for all $t$.

Consider the martingale
\begin{math}
  M\t := \sum_{s=0}^{t-1} \left(X^{s+1} - \Ext[s]{ \X^{s+1} }\right),
\end{math}
and the previsible process
\begin{math}
  A\t := \sum_{s=0}^{t-1} \left(\Ext[s]{ \X^{s+1} } - \X^s\right).
\end{math}
The following identity
\begin{equation}
  \X\t = \X^0 + A\t + M\t \label{X=AM} 
\end{equation}
is called the Doob decomposition of $\X\t$.  We refer to the terms of
$A\t$ as \emph{drifts} of $\X\t$ and to the terms in $M\t$ as
\emph{deviations}.

On account of the bound $X\t \leq 3X\t[0]$, we obtain from \qr{driftt}
that
\begin{math}
  A\t \leq a \, 3 \X^0.
\end{math}
and, from \qr{devv}, that
\begin{math}
  |M\tr - \Ext{ M\tr }| \leq \xi \, 3\X^0.
\end{math}
Since $M\t$ is a martingale, the Azuma-Hoeffding inequality implies
that
\begin{equation*}
  \Pr{|M\t| > \lambda} < 
  \exp\left\{ - \lambda^2 / 4 (\xi\sqrt{n} \, 3 \X^0)^2\right\}
\end{equation*}
if
\begin{equation*}
  \lambda = \OrdoTheta{ \xi\sqrt{n}(\log n)}\cdot\X\t[0]. 
\end{equation*}
It is then also seen that $|M\t| > \lambda$ is an event of probability
of order $n^{-\ordoomega{1}}$.

Hence, the probability $\X\t - X\t[0] > \Phi \X\t[0]$ where
\begin{displaymath}
  \Phi \X\t[0] = 
  \OrdoOmega{a\log n} \X\t[0] + \OrdoOmega{\xi\sqrt{n}\log n} \X^0 
\end{displaymath}
is of order $n^{-\ordoomega{1}}$; in the display above the first term
is greater than $A^t \leq 3a\X\t[0]$ and the second term bounds $\M\t$
within very high probability.  This proves \qr{uub}.  The stability of
\qr{lub} follows in the same manner from the stability of the event
$M\t > - \lambda$.  (In this case $A\t = 0$, since $\X$ is a
martingale.)
\end{proof}

\subsubsection{An estimate of $\xi$ for a certain type of sums}

\newcommand{\xmax}{\mathfrak{m}}
\newcommand{\killmax}{\mathcal{N}}

The inequalites we are dealing with have a common form and we will
repeatedly use the formula in \qr{cestimate} below in order to to estimate
the deviation parameter $\xi$ used in \lemref{concentration}.

The processes we consider are sums 
\begin{displaymath}
        \X\t = \sum_{i\in \mathcal{J}} \X_i\t, \qquad \X_{i}\t \geq 0
\end{displaymath}
for some index set $\mathcal{J}$. The terms have uniform bounds 
\begin{displaymath}
        0 \leq \X_i\t \leq \xmax,
\end{displaymath}
for some asymptotic expression $\xmax$.  It will also be the case that each
term $\X_i\t \geq 0$ in the sum changes moderately in the following sense:
For $\xi_{+}, \xi_{-} > 0$, we have for all $t\in\Time$ and
$i\in\mathcal{J}$ that
\begin{align}
                (1-\xi_-) \X_i\t &\leq \X_i\tr \leq (1+\xi_+) \X_i\t,
\end{align}
unless the term $\X_i\t$ is \emph{killed}, \ie unless $\X_i\tr = 0$
but $\X_i\t>0$.  We assume that the maximum number of such terms
killed is furthermore given by $\killmax$.

The following lemma is then immediate.
\begin{lem}\label{lemces}
  With $\X$, $\xmax$, $\killmax$, $\xi_{+}$ and $\xi_{-}$ as above
  the conclusions of \lemref{concentration} hold with
\begin{equation}\label{cestimate}
        \xi =  \frac{\killmax\,\xmax}{X^0} + \xi_- + \xi_+.
\end{equation} 
\end{lem}

\subsection{The stability of \qr{norm1}}\label{stabnorm}

Let $\ell(x)$ be any local line, \ie $\ell(x)=\ell_\RK(x)$ or
$\ell(x)=\ell_\RS(x)$.  For $t\in \Time$, let
\begin{equation}\label{Xnorm}
        X\t:= \sum_{y \in \ell(x)} \p\t(y).   
\end{equation} 
Our objective is to show that $(1-\Phi)\leq X^t \leq (1+\Phi)$ with
very high probability, where $\Phi = \log n/\sqrt n$.  By the
martingale property \qr{martingale}, the drift, $a$, for $\X^s$ is
zero and we have that $\X^0 = 1$.  Unless $\p\tr(y)=0$ we have, by
\qr{pincr}, that $\p\t(y) \leq \p\tr(y) \leq (1 + \pmax)\, \p\t(y)$.
Thus, in the notation from \lemref{lemces}, we have $\xi_- = 0$ and
$\xi_+ = \pmax$.  Moreover, each term in \qr{Xnorm} is smaller than
$\pmax$ (by \qr{pbnd}), so we have $\xmax = \pmax$ and hence, by
\lemref{lemces}, $\xi = \killmax\,\pmax +\ 0 + \pmax$.

Thus, in order to show that $\xi=\pmax$, it only remains to show that
$\killmax$ --- the maximum number of terms
``killed'' --- is of order $\Ordo{1}$. We claim that $\killmax \leq
2$. The number of terms killed is $\sum_{y\in\ell(x)} \K\tr(y)$ where
$\K\tr = \L\tr\circ\rok\t + \L\tr\circ\rod\t$. Moreover, the maps
$y\to z=\rok\t(y)$ and $y\to z=\rod\t(y)$ maps $y\in\ell(x)$
one-to-one into a corresponding local line $z\in\ell(\rok\t(y))$ and
$z\in\ell(\rod\t(y))$.  Since $\L\tr$ is latin this means that
\begin{equation}\label{uuux}
    \killmax \leq 
    \sum_{y} \L\tr\circ\rok\t(y)  +
    \sum_{y} \L\tr\circ\rod\t(y)  \leq 1 + 1  = 2.
\end{equation}
\qed

\subsection{Proof of stability for \qr{qrnd1}}

For a fixed $i \in \Rows$ and $k, l \in\Cols$, $k\not= l$, we consider the sum
$\X\t:= \sum_\g \X_\g\t$ where for $\g\in\Syms$
\begin{equation}
        \X_\g\t := \p\t(ik\g)\cdot \p\t(il\g).
\end{equation}
We have $\X^0 = 1/n$ and $\X_\g \leq \pmax^2$ by \qr{pbnd}. Our
objective is to prove that with very high probability $X\t \leq
(1+\Phi)X\t[0]$ where $\Phi = \log n /\sqrt n$. By,
\lemref{concentration}, it is enough to show that $\xi=\pmax$ and
$a=\pmax$. These bounds are proved in \qr{devXbnd} and \qr{aestim}
below.

By \qr{pincr}, each term $\X_\g\tr$ will either increase with a factor
at most $(1 + \pmax)^2 = 1 + \pmax$ or be killed.  Furthermore, at
most four terms can be killed, since, by the computation already done
in \qr{uuux}, at most two points in each of the cells $(i,k)$ and
$(i,l)$ are killed.  Hence, with $\xi$ as in \lemref{lemces}
\begin{equation}\label{devXbnd}
        \xi \leq \frac{4\pmax^2}{1/n} + 0 + \pmax = \pmax.
\end{equation}

In order to calculate the drift $a$ of $X\t$, fix $\g\in\Syms$ and let
$K_k = \K\tr(\ikg)$ and $r_k = \Ext{K_k}$. Note that $r_k =
\q\t\circ\rok(\ikg) + \q\t\circ\rod(\ikg) = \pmax$.  Then
\begin{equation}\begin{split}
  \Ext{ X_{\g}\tr}  = X_\g\t \cdot 
        \frac{\Ext{(1-K_k)(1-K_l)}}{ \Ext{(1-K_k)}\Ext{(1-K_l)} }  = \\
        =  X_\g\t \cdot  
        \frac{ 1 - r_k - r_l + \Ext{K_k K_l} }
        { 1 - r_k - r_l + r_k r_l }. \qquad
        \end{split}
        \label{driftX}
\end{equation} 
But, the event $K_k K_l \not= 0$ can happen only when the diagonal
through the cell $(i,k)$ intersects row $t+1$ in the column number
$l$ and vice versa and this can be the case for at most two values of
$t$.  For all other values of $t$ the diagonals and columns through
cells $(i,k)$ and $(i,l)$ intersect row $t+1$ at four disjoint
positions.  For these $t$, at most one cell in row $t+1$ is coloured by
any colour $\g$, so at most one of the indicators $K_k$ and $K_l$
equals one making $K_k K_l \equiv 0$.  Hence it holds, for all but at
most two values of $t$ and \emph{uniformly for all $\g$}, that
$\Ext{K_k K_l} = 0$ and from \qr{driftX} it is clear that $\Ext{
  X_{\g}\tr} \leq X_\g\t$ for these $t$.

Also, for the two possible exceptional values of $t$, when $\Ext{K_k
  K_l}$ is positive, we clearly have that
\begin{displaymath}
        \Ext{K_k K_l} \leq \min\{r_k, r_l\} \leq \pmax,
\end{displaymath}
even assuming the worst possible correlation. From \qr{driftX}, we
deduce that $\Ext{ X\tr} \leq X\t + \alpha \max\{X^0,X\t\}$, with
$\alpha=\pmax$.  Putting this together gives that $a = \sum_{s=0}^{t-1}
\alpha\t[s]$ as in \lemref{concentration} is bounded by
\begin{equation}\label{aestim}
        a \leq (t-2) \cdot 0 + 2 \pmax = \pmax.
\end{equation}
\qed

\subsection{Proving that \qr{pbnd1} is stable}

\newcommand{\so}[2][\ell]{{#1\cap\{\ctime <#2\}}}
\newcommand{\sumso}[2][\ell]{\sum_{\so[#1]{#2}}}
\newcommand{\Rtp}[1][\ell]{S} 
\newcommand{\that}{{{{\hat t}}}}
\renewcommand{\so}[1][t_0]{\ell^{<t_0}}

The derivation of the stability of \qr{pbnd1} is a bit more involved.
The important part is perhaps the identity \qr{iidd} below, which
makes it possible to relate the growth of the product $\prod_{s =
  0}^{t} (1 - \p^s\circ\ro^s)^{-1}$ to that of a sum $\Rtp\t =
\sum_{s=0}^{t-1} \tp\t[s]\circ\ro^s$ with suitable concentration
properties. Our aim is to prove that with very high probability
\begin{equation}\label{goal00}
  \frac {\p^t(x)}{\p^0(x)} \leq 
  (1+\ordo 1) \, \left(1 - (1+\ordo 1)\,\frac tn\right)^{-2}. 
\end{equation}
Since $t/n \leq (1-\e) = 1-\OrdoOmega 1$, the inequality \qr{goal00}
implies \qr{pbnd1} for large enough $n$.

Let $\ell$ denote a central parallell class, \ie $\ell$ is either of
$\ell_\KS$ or $\ell_\DS$.  Define the vector of stopping times $\Tc =
\Tc(x)$, $x\in\RKS$, giving the time that the central line $\ell(x)$
is ``killed'', \ie let
\begin{equation}\label{Tcdef}
  \Tc(x) := \inf\Setof{t \st \sum_{y\in\ell(x)} \L^{t}(y) = 1}
  \cup\{\infty\}.  
\end{equation}
Also, let
\begin{displaymath}
  \that(x,t) := t \wedge  (\Tc[\ell_\KS](x) \wedge \Tc[\ell_\DS](x)-1)
  \wedge \T \wedge\ctime(x). 
\end{displaymath}
Note that the value of $\that(x,t)$ is determined at time $t$ for all
$t$ and that $\p^\that(x)=\p^{\that(x,t)}(x)$ is an adapted process
which is increasing in $t$.  Moreover, we have $\p\t(x) = \p^\that(x)$
unless $t \geq \Tc[\ell_\KS](x) \wedge \Tc[\ell_\DS](x)$ and
$\Tc[\ell_\KS](x) \wedge \Tc[\ell_\DS](x) < \T \wedge\ctime(x)$, in
which case $\p\t(x) = 0$.

From the definition \qr{pdef}, we know that
\begin{multline}\label{expand00}
  \sup_{s\leq t} \frac {\p^s}{\p^0} = \prod_{s = 0}^{\that - 1}
  \frac 1{1 - \q^s\circ\rok^s - \q^s\circ\rod^s} =  
  \prod_{s = 0}^{\that - 1} \frac{(1 - \p^s\circ\rok^s)(1 -
    \p^s\circ\rod^s)} {1 - \q^s\circ\rok^s - \q^s\circ\rod^s} \, \cdot
  \\
  \prod_{\ell=\ell_\KS,\ell_\DS} \prod_{s = 0}^{\that - 1} (1 - \p^s\circ\ro^s)^{-1}.
\end{multline}
Note that the first factor to the right is negligible: By \qr{pbnd}
and the estimate in \qr{exstr}, it is of order $(1 + \qpd\pmax +
\pmax^2)^\that = (1+\qpd)$.  The aim will thus be to show that
\begin{equation}\label{goal001}
  \prod_{s = 0}^{\that - 1} (1 - \p^s\circ\ro^s)^{-1} 
  = (1- (1+\ordo1)\, \frac{\that}n)^{-1}.
\end{equation}
where $\ell=\ell_{\KS}$ or $\ell=\ell_{\DS}$.

For $\ell=\ell_{\KS}$ or $\ell=\ell_{\DS}$, define recursively the
adapted process $\tp\t \in \UI^\RKS$, $t\in \Time$, as follows: Let
$\tp^0(x) = \p^0(x) = \frac 1n$ and set
\begin{equation}\label{deftp}
  \tp\tr(x) :=
  \begin{cases}
    (1 - \Rtp\tr(x))\cdot\p\tr(x)  & t \leq (\Tc(x)-1)\wedge \T\wedge
    \ctime(x)  \\
    \tp\t(x) & \text{otherwise}
  \end{cases}
\end{equation}
where
$$\Rtp\t(x) := \sum_{s=0}^{t-1} \tp^s\circ\ro^s(x). $$  
The definition of $\tp\t$ in \qr{deftp} implies, for all $y \in
\ell(x)$, that either $\tp\tr(y)=0$, $\tp\tr(y)=\tp\t(y)$ or else
$$ \tp\tr(y) =  
\frac {(1-\p\t\circ\ro^s(y))\,\tp\t(y)}{1 - \q^s\circ\rok^s(y) -
  \q^s\circ\rod^s(y)}. $$ 
Therefore,
\begin{equation}\label{devXpb}
         \tp\t(y)
         \cdot (1-\pmax) \leq \tp\tr(y) < \tp\t(y) \cdot (1+\pmax), 
\end{equation}
except for the case when the term $\tp\t(y)$ is killed.  Note also
that killing of $\tp\t(y)$ occurs exactly when $\Tc[\ell'](y) = t+1 <
\Tc[\ell](y)$, where $\ell'(y)$ is the \emph{complementary} central line, \ie
$\ell' = \ell_{\KS}$ if $\ell = \ell_{\DS}$ and vice versa.

If $t=\that$ then $\p\t = \tp\t/(1-\Rtp\t)$ and we deduce the identity
\begin{equation}\label{iidd}
    \prod_{s = 0}^{\that - 1} (1 - \p^s\circ\ro^s)^{-1} =
    \prod_{s = 0}^{\that - 1} \frac{1-\Rtp^{s}}{1-\Rtp^{s+1}} =
    \frac 1{1-\Rtp^\that}.
\end{equation}
In order to prove the stability of \qr{goal00} it is therefore enough
to show that for all $x$ with very high probability
\begin{equation}
    \Rtp\t[t_0](x) \leq (1+\ordo{1}) \,(t_0/n),
    \label{goal0}
\end{equation}
where the time $t_0\in\Time$ is fixed but arbitrary. 

So fix $t_0\in\Time$ and let $x\in\RKS$ be arbitrary. In order to ease
the notation, we from now on suppress the dependency of $x\in\RKS$ for
most quantitites.  Let $\so(x) := \Setof{y\in\ell(x)\st \ctime(y)<
  t_0}$.  Define for $t\in [0,t]$ the variable $X\t\in \RR_+^{\RKS}$
by
\begin{displaymath}
        X\t(x) := \sum_{y\in \so(x)} \tp\t(y).   
\end{displaymath} 
Since $\tp\t = \tp^{t\wedge \ctime}$, we have that $X^{t_0} =
\Rtp^{t_0}$. 
Since
$|\so| = t_0$ and $\tp^0 = 1/n$ we have $X^0 = t_0/n$. Moreover, on
account of \qr{devXpb}, we have $X\t \leq (1+ t\pmax)X^0$ and thus the
sought bound, \qr{goal0}, clearly follows, if $t_0 = \ordo{n}$. We may
therefore assume that $t_0$ is greater than, say, $n^{2/3}$ and hence
that $X^0 \geq n^{-1/3}$.

Thus, we can conclude the proof by showing that the total drift, $a$,
of $\X\t$ satisfies $a \leq \qpd$ and that the deviation parameter,
$\xi$, satisfies $\xi=\Ordo{n^{-2/3}}$. Then, from
\lemref{concentration}, we can deduce that $$
\Rtp^{t_0}=\X\t[t_0] \leq \left(1+\Ordo{n^{-1/6}\,\log
    n}\right)\,(t_0/n), $$
with very high probability. 

\newcommand{\roc}{{\ell'}}
From the definition \qr{deftp}, we know that $\Ext{\tp\tr/\tp\t}$ equals $1$
in the case when $t \geq \Tc\wedge\T\wedge\ctime$ and that otherwise it equals 
\begin{multline*}
        \Prt{ \Tc = t+1 } \cdot 1 + \Prt{ \Tc > t+1 } \cdot 
        \frac {1 - \Rtp\tr}{1-\Rtp\t} 
        \cdot \Ext{{\p\tr}/{\p\t} \mid \Tc > t+1} \\
        = \q\t\circ\ro\t + (1 - \q\t\circ\ro\t) \cdot \frac 
        {1 - \p\t \circ\ro\t}{1 - \q\t\circ\ro\t} 
        = 1 + \q\t\circ\ro\t - \p\t\circ\ro\t, 
\end{multline*} 
where we have used that $\frac {1 - \Rtp\tr}{1-\Rtp\t} = {1 - \p\t 
\circ\ro\t}$ and that $\Ext{{\p\tr}/{\p\t} \mid \Tc > t+1}$ equals 
\begin{equation*}
\frac{\q\t\circ\roc\t}{1-\q\t\circ\ro\t} \cdot 0 + 
     (1-\frac{\q\t\circ\roc\t}{1-\q\t\circ\ro\t}) 
     \cdot \frac 1{1 - \q\t\circ\ro\t-\q\t\circ\roc\t} =
     \frac1{1-\q\t\circ\ro\t},
\end{equation*} 
where $\ell'$ is the complementary line. 

Hence, for all $t$, 
\begin{equation}\label{tpex}
    \begin{aligned}
        \Ext{\tp\tr} & = 
        \begin{cases}
        \tp\t \cdot (1 + \q\t\circ\ro\t - \p\t\circ\ro\t) 
                  & t < \ctime \wedge \T \wedge \Tc \\
        \tp\t & \text{otherwise}
        \end{cases}\\
        &\leq \tp\t \cdot (1 + \pmax \qpd).
    \end{aligned}
\end{equation}
From \qr{tpex} we can estimate the drift term 
\begin{equation*}\label{aesta}
        a \leq t_0 \qpd\pmax = \qpd. 
\end{equation*}

Note that only one cell in the row $t+1$ is coloured with the
colour $\g$ common to all points in $\ell(x)$ and that this cell can
lie on at most one crossing central line $\ell'(y)$ intersecting the
given central line $\ell(x)$.  This means that only one term in the
sum $\sum_{\so} \tp\t$ above can be killed at a time.  Putting this
and the estimate \qr{devXpb} into the formula \qr{cestimate} gives
\begin{equation*}\label{xieste}
        \xi \leq  1\cdot \pmax/X^0 + \pmax + \pmax  = \Ordo{ n^{-2/3} }
\end{equation*}
\qed


\ifx\undefined\bysame
\newcommand{\bysame}{\leavevmode\hbox to3em{\hrulefill}\,}
\fi

\end{document}